# MONOTONE PROPERTIES OF RANDOM GEOMETRIC GRAPHS HAVE SHARP THRESHOLDS

By Ashish Goel,[1] Sanatan Rai[2] and Bhaskar Krishnamachari

*Stanford University, Stanford University and University of Southern California*

Random geometric graphs result from taking $n$ uniformly distributed points in the unit cube, $[0,1]^d$, and connecting two points if their Euclidean distance is at most $r$, for some prescribed $r$. We show that monotone properties for this class of graphs have sharp thresholds by reducing the problem to bounding the bottleneck matching on two sets of $n$ points distributed uniformly in $[0,1]^d$. We present upper bounds on the threshold width, and show that our bound is sharp for $d=1$ and at most a sublogarithmic factor away for $d \geq 2$. Interestingly, the threshold width is much sharper for random geometric graphs than for Bernoulli random graphs. Further, a random geometric graph is shown to be a subgraph, with high probability, of another independently drawn random geometric graph with a slightly larger radius; this property is shown to have no analogue for Bernoulli random graphs.

**1. Introduction.** Consider $n$ points distributed uniformly and independently in the unit cube $[0,1]^d$. Given a fixed distance $r > 0$, connect two points if their Euclidean distance is at most $r$. Such graphs are called random geometric graphs, and are denoted by $G^{(d)}(\mathcal{X}_n; r)$, as in [24]. Classically, these graphs have been the subject of much study because of connections to percolation, statistical physics, hypothesis testing and cluster analysis. Further, random geometric graphs are better suited than more combinatorial classes (such as Bernoulli random graphs) to model problems where the existence of an edge between two different nodes depends on their spatial

Received April 2004; revised May 2005.
[1]Supported in part by NSF CAREER Award CCR-0339262 and by an Alfred P. Sloan faculty fellowship.
[2]Supported by NSF Award CCR-0339262.

distance. As a result, random geometric graphs have received increased attention in recent years in the context of distributed wireless networks (such as sensor networks), see, for example, [12, 13, 14]; and layout problems as in [5, 16, 24]. Another area is cluster analysis, especially its applications in medicine, biology and ecology; these may be found in [10].

In applications such as distributed wireless networks, the connectivity of random geometric graphs is of interest. Gupta and Kumar showed that for $d = 2$, if $\pi r(n)^2 = (\log n + c_n)/n$, then as $n \uparrow \infty$ the graph is connected almost surely as $n \uparrow \infty$ if $c_n \uparrow \infty$ and is disconnected almost surely if $c_n \downarrow -\infty$ [12]. This result is remarkably similar to the corresponding result for Bernoulli random graphs (also known as Erdős–Renyi graphs). An instance of a Bernoulli random graph is obtained by taking $n$ points and connecting any two with probability $p$, independently of all other pairs. This class of graphs is denoted by $\mathcal{G}_{n,p}$. Erdős and Renyi [7, 8] showed that if $p(n) = (\log n + c_n)/n$, then the graph is a.s. connected or disconnected as $c_n \uparrow \infty$ or $c_n \downarrow -\infty$. For $d = 2$, Gupta and Kumar's result can also be obtained from Penrose's work on the longest edge of minimal spanning tree of $G(\mathcal{X}_n; r(n))$ [23]. Connectivity of random geometric graphs for $d = 1$ was also studied by Godehardt and Jaworski [11]. Connectivity results under the $l_\infty$-norm may be found in [1].

In both random geometric graphs and Bernoulli random graphs, property thresholds are of great interest. To quote Bollobás [2]:

> One of the main aims of the theory of random graphs is to determine when a given property is likely to appear.

Particularly interesting are thresholds for monotone properties, of which connectivity is a classic example. A seminal result of Friedgut and Kalai [9] states that all monotone graph properties have a sharp threshold in Bernoulli random graphs, and the threshold width is $\delta(\varepsilon) = \mathcal{O}(\log \varepsilon^{-1}/\log n)$. They also demonstrated a monotone property with a threshold width of $\Omega(1/\log^2 n)$ and conjectured that this is tight [i.e., the best upper bound on threshold width is $\mathcal{O}(1/\log^2 n)$]. Their upper bound on threshold width was improved to $\mathcal{O}(1/\log^{2-\gamma} n)$ for all $\gamma > 0$ by Bourgain and Kalai [4].

The similarity of the connectivity threshold for random geometric graphs and Bernoulli random graphs led to the conjecture that all monotone properties also have a sharp threshold in random geometric graphs (see [18, 20] for a more detailed discussion). For the $d = 1$ case, sharp thresholds for monotone properties are implicit in the recent work of McColm [20], though he does not compute bounds on the width. The definition of sharp thresholds we use in this paper is the one used by Friedgut and Kalai and is based on the threshold width. The definition used by McColm is the one used in the text by Janson, Łuczak and Ruciński [17], and is stronger than the one used by Friedgut and Kalai; we discuss this in more detail at the end of the Introduction. The analysis of random geometric graphs is technically challenging because of dependence of the edges. The triangle inequality implies that the event that points $x$ and $z$ are connected is not independent of the event $\{(x,y)$ and $(y,z)$ are edges$\}$. This is in stark contrast to the case of Bernoulli random graphs. Hence proof techniques that have been successful for $\mathcal{G}_{n,p}$ cannot be exploited in the case of random geometric graphs.

*Our results.* We show that all monotone graph properties have a sharp threshold for random geometric graphs, thus resolving the above conjecture. In fact, the threshold width for random geometric graphs is much sharper than for Bernoulli random graphs. In order to state our results formally, we need to establish some notation and some definitions.

We use the symbol $\sim$ to mean "distributed as," so that $G \sim G^{(d)}(\mathcal{X}_n; r)$ means that $G$ is picked from $G^{(d)}(\mathcal{X}_n; r)$. For ease of notation, we omit the superscript $d$ in $G^{(d)}(\mathcal{X}_n; r)$ as the dimension will be clear from the context. The critical radius for connectivity is defined as $r_c := (\log n / \pi_d n)^{1/d}$, where $\pi_d$ is the volume of the unit sphere in $\mathbb{R}^d$.

A graph property $A$ is a set of undirected and unlabelled graphs. A property $A$ is *increasing* if and only if

$$G \in A \implies (\forall G')[(V(G') = V(G) \text{ and } E(G) \subseteq E(G')) \Rightarrow G' \in A].$$

Intuitively speaking, an increasing property is one that is preserved when edges are added to the graph. A graph property $A$ is *monotone* if either $A$ or $A^c$ is increasing. Without loss of generality, for the rest of the paper, we shall implicitly mean increasing properties when referring to monotone properties.

If $A$ is an increasing property, then for $0 < \varepsilon < 1/2$, let $r(n,\varepsilon) = \inf\{r > 0 : \mathbb{P}\{G(\mathcal{X}_n; r) \in A\} \geq \varepsilon\}$. Define further $\delta(n,\varepsilon) := r(n, 1-\varepsilon) - r(n,\varepsilon)$. A property is said to have a *sharp threshold* if $\delta(n,\varepsilon) = o(1)$ for all $0 < \varepsilon < 1/2$.

Our main results are:

THEOREM 1.1. *For every monotone property, the width $\delta(n,\varepsilon)$ is*

$$\mathcal{O}\left(\sqrt{\frac{\log \varepsilon^{-1}}{n}}\right)$$

*for $d = 1$. For $d = 2$,*

$$\delta(n,\varepsilon) = \mathcal{O}(r_c \log^{1/4} n) \equiv \mathcal{O}(\log^{3/4} n / \sqrt{n}),$$

*and for $d \geq 3$, the width*

$$\delta(n,\varepsilon) = \mathcal{O}(r_c) \equiv \mathcal{O}(\log^{1/d} n / n^{1/d}).$$

Thus, all monotone properties have sharp thresholds. Observe that the width is much sharper than the threshold width for $\mathcal{G}_{n,p}$. Moreover, we prove a stronger result: the graphs $G(\mathcal{X}_n; r)$ become subgraphs of $G(\mathcal{X}_n; \rho)$ for $\rho > r$, with hight probability (w.h.p.) as $n \uparrow \infty$. Note only that this is not the case for Bernoulli random graphs—we shall make this precise in Section 2.

For the lower bounds we have:

THEOREM 1.2. *For $d \geq 2$, there exists a monotone property with width $\delta(n, \varepsilon) = \Omega((\log \varepsilon^{-1})^{1/d} n^{-1/d})$. For $d = 1$, there exists a monotone property with width*

$$\delta(n, \varepsilon) = \Omega(\sqrt{\log \varepsilon^{-1}}/\sqrt{n}).$$

Hence, we have a tight characterization of the threshold width for $d = 1$, and our upper bounds are only a sublogarithmic factor away for $d \geq 2$.

The key idea is to relate the behavior of monotone properties to the weight of the "bottleneck" matching (to be defined later) of the bipartite graph whose vertex sets are obtained by distributing $n$ points uniformly and independently in $[0,1]^d$. Sharp results on the "bottleneck" matching weight are implicit in the work of Leighton and Shor [19] for $d = 2$ and Shor and Yukich [26] for $d \geq 3$. We repeat them here for convenience.

THEOREM 1.3 ([19, 26]). *Consider the bipartite graph on $2n$ points, where each set of $n$ points is distributed uniformly and independently in the unit cube $[0,1]^d$, and independently of each other. If $M$ is the length of the bottleneck matching, then w.h.p. as $n \uparrow \infty$:*

$$M = \begin{cases} \Theta(r_c \log^{1/4} n), & \text{if } d = 2, \\ \Theta(r_c), & \text{if } d \geq 3. \end{cases}$$

In Section 4 we present our own proof of the bound for $d \geq 3$. The proof for the $d \geq 3$ case in Shor and Yukich [26] invokes results from polyhedral geometry. We present a simpler proof that relies only on the properties of order statistics and Chernoff bounds. We prove that for $d = 1$, the bottleneck matching is $\mathcal{O}(\sqrt{\log \varepsilon^{-1}}/\sqrt{n})$ with probability $1 - \varepsilon$. Moreover, our results also hold for any $\ell_p$-norm when $p > 1$, and not just under the Euclidean norm as in the setting of this paper. We omit the details, as they require straightforward modifications of the proofs given herein.

It might seem curious that we do not report the dependence on $\varepsilon$ in some of our bounds on the threshold width. This is because the results of Shor and Yukich as well as those of Leighton and Shor are high probability results: the bottleneck matching length is $\Theta(r_c \log^{1/4} n)$ in two dimensions and $\Theta(r_c)$ in higher dimensions not just in expectation but with probability $1 - o(n^{-\beta})$ for some $\beta > 0$. Hence, in asymptotic notation, our upper bound on $\delta(n, \varepsilon)$ in two and higher dimensions does not depend on $\varepsilon$ as long as $\varepsilon = \Omega(n^{-\beta})$.

*Related work.* There is a vast body of literature that is directly related to this paper. It would require a survey paper to even mention the salient results with any degree of honesty. We can only point the reader to the book by Penrose [24], the papers by Gupta and Kumar [12, 13, 14] and the paper by Shakkottai, Srikant and Shroff [25]. We note here that our techniques imply a sharp threshold for the coverage problem as discussed in [25], which is not a graph problem. We omit the details. The theory of Bernoulli random graphs is covered in the books by Bollobás [2] and Janson, Łuzak and Rucinski [17]. For some results on matchings in a similar context, see the paper by Holroyd and Peres [15] and for some results on covering algorithms see [3]. Sharp thresholds for random geometric graphs were conjectured in [18] and [20]. Muthukrishnan and Pandurangan [21] obtained asymptotically tight thresholds for connectivity, covering and routing-stretch in $d$ dimensions using a new technique called *bin-covering*.

*Additive versus multiplicative thresholds.* In this paper we are primarily concerned with bounding the threshold width of a property, along the lines of Friedgut and Kalai [9]. Informally, this corresponds to proving sharp "additive" thresholds. As we mentioned earlier, the notion of sharp thresholds presented in [17] or in [20] is stronger in that they require $\delta(n,\varepsilon)/r_\Pi(n) = o(1)$. Informally, this corresponds to "multiplicative" thresholds. We observe that our Theorem 1.1 also yields sharp thresholds in this stronger sense, provided the threshold radius is high enough. More precisely, if $\Pi$ is a monotone property and $r_\Pi(n)$ is its threshold radius, such that:

$$\begin{aligned} r_c &= o(r_\Pi(n)) & \text{when } d \geq 3, \\ r_c &= o(r_\Pi(n)/\log^{1/4} n) & \text{when } d = 2, \\ \sqrt{n}\, r_\Pi(n) &\to \infty & \text{when } d = 1, \end{aligned}$$

then $\Pi$ also has a sharp threshold in the sense of Janson, Łuczak andRuciński [17].

*Plan of this paper.* We first establish the relationship between monotone properties and bottleneck matchings and prove the upper bound in Section 2. In Section 3 we furnish the lower bounds. In Section 4 we discuss the upper bound for $d \geq 3$, and in Section 5 for $d = 1$. We conclude in Section 6 with some open problems.

**2. Bottleneck matchings and monotone properties.** Recall that in a bipartite graph with vertex sets $V_1$ and $V_2$, a perfect matching is a bijection $\phi: V_1 \to V_2$, such that each $v$ is adjacent to $\phi(v)$. Thus a perfect matching is a disjoint collection of edges that covers every vertex. If the graph is weighted, then we define the *weight* of the matching as the maximum weight of any edge in the matching. A *bottleneck* matching is the perfect matching with the minimum weight.

Let $S_1$ and $S_2$ denote two sets of $n$ points each, where the points are i.i.d., chosen uniformly at random from the set $[0,1]^d$. Form the complete bipartite graph on $(S_1, S_2)$ and let the weight of an edge be the Euclidean distance between its endpoints. Let $M_n^{(d)}$ denote the bottleneck matching weight of this graph. We omit the dimension $d$ where it is clear from the context.

We first link the weight of the bottleneck matching with a containment property on random geometric graphs. We shall write $G \subset G'$ to mean that the graph $G$ is contained in the graph $G'$, that is, is isomorphic to a subgraph of $G'$.

LEMMA 2.1. *Suppose $\mathbb{P}\{M_n > \gamma(n)\} \leq p$ for some function $\gamma(n)$ and some constant $p$. For any radius $r$, consider independent random samples $G \sim G(\mathcal{X}_n; r)$ and $G' \sim G(\mathcal{X}_n; r + 2\gamma(n))$ in $d$ dimensions. Then, $\mathbb{P}\{G \subset G'\} \geq 1 - p$.*

PROOF. Let $V$ represent the set of points in graph $G$ and $V'$ the set of points in graph $G'$. Let $\phi$ denote the bottleneck matching between $V$ and $V'$; then $M_n$ is the weight of this matching. Suppose $(u,v) \in E(G)$, that is, $\|u-v\|_2 \leq r$. Then, using triangle inequality,

$$\|\phi(u) - \phi(v)\|_2 \leq \|\phi(u) - u\|_2 + \|u - v\|_2 + \|v - \phi(v)\|_2.$$

But $\|\phi(u) - u\|_2$ and $\|\phi(v) - v\|_2$ are both at most $M_n$, and hence $\|\phi(u) - \phi(v)\|_2 \leq 2M_n + r$. If $M_n \leq \gamma(n)$, then the mapping $\phi$ establishes that $G \subset G'$, and hence $\mathbb{P}\{G \subset G'\} \geq 1 - p$. □

The main result linking monotone properties to bottleneck matchings is:

THEOREM 2.2. *If $\mathbb{P}\{M_n > \gamma(n)\} \leq p$, then the $\sqrt{p}$-width of any monotone property in $d$ dimensions is at most $2\gamma(n)$.*

PROOF. Let $p = \varepsilon^2$, so that $\mathbb{P}\{M_n > \gamma(n)\} \leq \varepsilon^2$. Let $\Pi$ be an arbitrary increasing monotone property. Let $r_L = r(n, \varepsilon)$, $r_U = r_L + 2\gamma(n)$. Let $G \sim G(\mathcal{X}_n; r_L)$, and $G' \sim G(\mathcal{X}_n; r_U)$, and define $q := \mathbb{P}\{G' \notin \Pi\}$. Since $G$ is independent of $G'$, $\mathbb{P}\{G \in \Pi, G' \notin \Pi\} = \varepsilon \cdot q$. The monotonicity of $\Pi$ implies that if $G \in \Pi$ and $G' \notin \Pi$, then $G \not\subset G'$. This means that $\mathbb{P}\{G \not\subset G'\} \geq \varepsilon \cdot q$. By Lemma 2.1 above, $\mathbb{P}\{G \not\subset G'\} \leq p$, so that we must have $q \leq \varepsilon$. But then $r(n, 1 - \varepsilon) \leq r(n, 1 - q) = r_U$, so that $\delta(n, \varepsilon) \leq r_U - r_L = 2\gamma(n)$. □

With Theorem 2.2, the upper bound theorem follows with very little more work:

PROOF OF THEOREM 1.1. Leighton and Shor [19] show that
$$M_n^{(2)} = \Theta(r_c \log^{1/4} n),$$
with probability at least $1 - n^{-\kappa}$, for some $\kappa > 0$, and Shor and Yukich [26] show that
$$M_n^{(d)} = \Theta(r_c) \qquad \text{for } d \geq 3,$$
with probability at least $1 - n^{-\kappa'}$, for some $\kappa' > 0$. Hence Theorem 2.2 implies that $\delta(n, \varepsilon) = \mathcal{O}(r_c \log^{1/4} n)$ for $d = 2$ and $\delta(n, \varepsilon) = \mathcal{O}(r_c)$ for $d \geq 3$, for any constant $\varepsilon > 0$. In fact, the bound on $\delta(n, \varepsilon)$ holds for any $\varepsilon = \Omega(n^{-c})$, where $c > 0$ is a constant.

In Proposition 5.1 (see Section 5), we show that for $d = 1$,
$$\mathbb{P}\left\{M_n^{(1)} \leq \frac{\beta}{\sqrt{n}}\right\} \geq 1 - \exp(-c\beta^2),$$
for some $c > 0$. By applying Theorem 2.2 with $\varepsilon = \exp(-c\beta^2)$ we obtain
$$\delta(n, \varepsilon) = \mathcal{O}\left(\sqrt{\frac{\log \varepsilon^{-1}}{n}}\right). \qquad \square$$

REMARK 1. We have in fact shown that $G(\mathcal{X}_n; r)$ is a subset of $G(\mathcal{X}_n; r')$ w.h.p., when $r' = r + o(1)$. The corresponding result does not hold for Bernoulli random graphs. To see this suppose that $G \sim \mathcal{G}_{n,p}$ and $G' \sim \mathcal{G}_{n,P}$. For $G \subset G'$, every edge in $G$ must exist in $G'$. Hence, for $M = \binom{n}{2}$, $q = 1 - p$ and $Q = 1 - P$, and a given matching $\phi$:

$$\mathbb{P}\{G \subset G' \text{ under } \phi\} = \sum_{K=0}^{M} \binom{M}{K} p^K q^{M-K}$$
$$\times \sum_{L=0}^{M-K} \binom{M-K}{L} P^{K+L} Q^{M-K-L}$$
$$= \sum_{K=0}^{M} \binom{M}{K} p^K q^{M-K}$$
$$\times P^K \sum_{L=0}^{M-K} \binom{M-K}{L} P^L Q^{M-K-L}$$
$$= \sum_{K=0}^{M} \binom{M}{K} (pP)^K q^{M-K} (P+Q)^{M-K}$$
$$= (pP + q)^M$$
$$= (p(1 - Q) + q)$$
$$= (1 - pQ)^M.$$

Choose $p = 1/4$, and $P = 3/4$. As there are $n!$ matchings:

$$\mathbb{P}\{G \subset G'\} \leq n! \exp\left(-\frac{n(n-1)}{32}\right).$$

The last expression goes to zero as $n \uparrow \infty$. Hence, even in this extreme case when $P - p = 1/2$, we do not have $\mathcal{G}_{n,p} \subset \mathcal{G}_{n,P}$ with constant probability.

**3. The lower bounds.** We now present examples of monotone properties to show that our bounds are tight in the $d = 1$ case and within a sublogarithmic factor for $d \geq 2$.

For the $d = 1$ case, we consider the property $\Pi$, defined by

$$G \in \Pi \iff \min_{u \in V} \deg(u) \geq \frac{|V|}{4},$$

where $V$ is the vertex set of $G$. Let $x_1, \ldots, x_n$ be the $n$ uniformly distributed points in $[0, 1]$; these are the vertices of the graph $G$. Let $x_{(i)}$ denote the $i$th order statistic. We have the following two estimates:

LEMMA 3.1. *If $0 < \varepsilon \leq 0.5 e^{-44/6}$, then for property $\Pi$:*

$$r(n, 1 - \varepsilon) \geq \frac{1}{4} + \sqrt{\frac{\log 1/2\varepsilon}{2n}}.$$

PROOF. Let $u = 1/4 + \Delta$, where $\Delta > 0$ is to be determined later. Pick $G \sim G(\mathcal{X}_n; u)$, and let the vertices be $x_1, \ldots, x_n$. Then $x_1, \ldots, x_n$ are distributed uniformly in $[0, 1]$. Observe that

$$\mathbb{P}\{x_{(n/4)} > u\} \geq \varepsilon \implies \mathbb{P}\left\{\deg(x_{(1)}) < \frac{n}{4}\right\} \geq \varepsilon$$
$$\implies \mathbb{P}\{G \notin \Pi\} \geq \varepsilon,$$

where in the first implication we have used the fact that $\deg(x_{(1)}) < n/4 \iff x_{(n/4+1)} - x_{(1)} > u$, and that $x_{(n/4+1)} - x_{(1)} \stackrel{d}{=} x_{(n/4)}$.

Now, $\mathbb{P}\{x_{(n/4)} > u\} = \mathbb{P}\{\mathsf{Bin}(n, u) < n/4\}$, and for some suitably large $n_0$,

$$\mathbb{P}\{\mathsf{Norm}(0, 1) < -\sqrt{n}\Delta\} \geq 2\varepsilon \implies \mathbb{P}\{\mathsf{Bin}(n, u) < n/4\} \geq \varepsilon,$$

whenever $n > n_0$, by the Normal approximation to the Binomial.

Put $\Delta = \beta/\sqrt{n}$ for $\beta = \sqrt{6 \log(0.5/\varepsilon)/11}$. Then for $0 < \varepsilon \leq 0.5 e^{-44/6}$, we have $\beta \geq 2$. With a little bit of work, one can see that $\beta^2/2 \leq -4\beta/3 - \log 2\varepsilon$. Since $x \geq \log x$ for $x \geq 1$, the last inequality implies that $\beta^2/2 \leq \log(3\beta^{-1}/4) - \log(2\varepsilon)$. Observe that any $\beta \geq 2$ satisfies

$$\frac{3}{4\beta} \leq \frac{1}{\beta} - \frac{1}{\beta^3}$$

so that
$$\frac{\beta^2}{2} \leq \log\left(\frac{1}{\beta} - \frac{1}{\beta^3}\right) + \log\frac{1}{2\varepsilon},$$
or
$$(\beta^{-1} - \beta^{-3})\exp\left(-\frac{\beta^2}{2}\right) \geq 2\varepsilon.$$

But then by Theorem 1.4 of [6], we can conclude $\mathbb{P}\{\mathsf{Norm}(0,1) > \beta\} \geq 2\varepsilon$. This shows that
$$\mathbb{P}\left\{G\left(\mathcal{X}_n; \frac{1}{4} + \sqrt{\frac{\log(2\varepsilon)^{-1}}{2n}}\right) \notin \Pi\right\} \geq \varepsilon,$$
and since $\Pi$ is increasing, this means that
$$r(n, 1-\varepsilon) \geq \frac{1}{4} + \sqrt{\frac{\log(2\varepsilon)^{-1}}{2n}}. \qquad \square$$

LEMMA 3.2. *For property $\Pi$:*
$$r\left(n, \frac{1}{2}\right) \leq \frac{1}{4} + \frac{c}{\sqrt{n}},$$
*for some fixed constant $c > 0$.*

PROOF. Suppose $G \sim G(\mathcal{X}_n, l)$. For any $u \in V(G)$, write $\deg_L(u)$ for the number of points to the left of $u$ and adjacent to it, and similarly let $\deg_R(u)$ stand for the number of right neighbors. Note that if $\deg(x_{(1)}) \geq n/4$, then necessarily $\deg(x_{(i)}) \geq n/4$ for all $1 \leq i \leq n/4$. With this observation we have:

$$\mathbb{P}\{G \notin \Pi\} = \mathbb{P}\left\{\bigcup_{1 \leq i \leq n}\left\{\deg(x_{(i)}) < \frac{n}{4}\right\}\right\}$$
$$\leq \mathbb{P}\left\{\left\{\deg(x_{(1)}) < \frac{n}{4} \text{ or } \deg(x_{(n)}) < \frac{n}{4}\right\}\right\}$$
$$+ \mathbb{P}\left\{\bigcup_{n/4 < i < 3n/4}\left\{\deg_L(x_{(i)}) < \frac{n}{8}\right\}\right\}$$
$$+ \mathbb{P}\left\{\bigcup_{n/4 < i < 3n/4}\left\{\deg_R(x_{(i)}) < \frac{n}{8}\right\}\right\}$$
$$\leq 2\underbrace{\mathbb{P}\left\{\deg(x_{(1)}) < \frac{n}{4}\right\}}_{(1)} + n\underbrace{\mathbb{P}\left\{\deg(x_{(i)}) < \frac{n}{8}\right\}}_{(2)}.$$

To bound the first term (1), first observe that by arguing as in the last lemma $\mathbb{P}\{\deg(x_{(1)}) < n/4\} = \mathbb{P}\{\mathsf{Bin}(n,l) < n/4\}$. By applying Chernoff bounds, we can find a $C_1 > 0$ so that $\mathbb{P}\{\mathsf{Bin}(n,l) < n/4\} < e^{-C_1^2/2}$, when $l = 1/4 + C_1/\sqrt{n}$.

To bound the second term (2), again apply Chernoff's bounds to find $C_2 > 0$, such that for $l = 1/4 + C_2/\sqrt{n}$,

$$n\mathbb{P}\left\{\deg(x_{(1)}) < \frac{n}{8}\right\} = n\mathbb{P}\left\{\mathsf{Bin}(n,l) < \frac{n}{8}\right\} \le n\exp\left(-\frac{n}{32}\right),$$

so that for $n$ large enough this term is overwhelmingly small. Therefore, for $c \ge \max(C_1, C_2)$, and $l = 1/4 + c/\sqrt{n}$, we have

$$\mathbb{P}\{G \notin \Pi\} \le e^{-c^2/2} + n\exp\left(-\frac{n}{32}\right),$$

so that $r(n, 1/2) \le 1/4 + c/\sqrt{n}$, for a suitably chosen $c$. □

Lemmas 3.1 and 3.2 show that for the graph property $\Pi$ defined by

$$G \in \Pi \iff \min_{u \in V} \deg(u) \ge \frac{|V|}{4},$$

we have

$$r_\Pi(n, 1-\varepsilon) \ge \frac{1}{4} + \sqrt{\frac{\log 1/2\varepsilon}{2n}} \quad \text{when } 0 < \varepsilon \le 0.5e^{-44/6},$$

$$r_\Pi\left(n, \frac{1}{2}\right) \le \frac{1}{4} + \frac{c}{\sqrt{n}}.$$

Hence we have shown the $d \ge 2$ case of Theorem 1.2:

PROOF OF THEOREM 1.2, LOWER BOUND FOR $d = 1$. Immediate from the last two lemmas. □

THEOREM 3.3. *For $d \ge 2$, there exists an increasing property $\Pi$ such that for $0 < \varepsilon < 1/2$, the threshold width satisfies*

$$\delta(n, \varepsilon) = \Omega(n^{-1/d}).$$

PROOF. Let $\Pi$ be the property that the graph is complete. This is trivially a monotone property.

Suppose $0 < \varepsilon < 1/2$. Let $u := \sqrt{d}(1 - 2\Delta_+)$ [see Figure 1(a)], for $\Delta_+$ such that

$$0 < \Delta_+ \le \frac{1}{n^{-1/d}}[\min(\sqrt{2\varepsilon}, \log 2/2)]^{1/d},$$

and pick $G \sim G(\mathcal{X}_n; u)$. Fix a pair of diagonally opposite corner cubes with side $\Delta_+$, and consider the event that there is exactly one point in each. If this happens, then the graph is not complete, since the points are more than $u$ apart. Hence:

$$\mathbb{P}\{G \notin \Pi\} \geq \binom{n}{2} (\Delta_+^d)^2 \cdot 2 \cdot (1 - 2\Delta_+^d)^{n-2},$$

since $\Delta_+ < (\log 2/(2n))^{1/d}$. Thus, for large enough $n$ we have

$$(1 - 2\Delta_+^d)^{n-2} \geq \left(1 - \frac{2\log 2}{2n}\right)^{n-2} \geq \frac{1}{2},$$

which implies that

$$\mathbb{P}\{G \notin \Pi\} \geq \binom{n}{2}(\Delta_+^d)^2 \geq \varepsilon.$$

The last inequality follows because we chose $\Delta_+^d \leq n^{-1} \times \min(\sqrt{2\varepsilon}, \log 2/2)$. Therefore,

$$(1) \qquad r(n, 1-\varepsilon) \geq u \geq \sqrt{d}\left(1 - \frac{c\varepsilon^{1/2d}}{n^{1/d}}\right).$$

Now we shall bound $r(n)$ above. To this end set $l = \sqrt{d - 1 + (1 - 4\Delta_-)^2}$ [see Figure 1(b)]. Now suppose that $\Delta_- = (\log \varepsilon^{-1})^{1/d}/(4n^{1/d})$, and pick $G \sim G(\mathcal{X}_n; l)$. Using elementary geometry it is easy to see that, if none of the $n$ points lies in any of the $2^d$ cubes of side $2\Delta_-$ at the corners of $[0,1]^d$, then the graph is complete. Hence, for $n$ large:

$$\mathbb{P}\{G \in \Pi\} \geq (1 - 2^d(2\Delta_-)^d)^n \geq \exp(-n(4\Delta_-)^d),$$

where we have used the fact that $1 - x \geq e^{-x}$, when $x \geq 0$, so that $(1-x)^n \geq e^{-nx}$, if $x < 1$. Since $\exp(-n(4\Delta_-)^d) = \varepsilon$, by our choice of $\Delta_-$, it must be that

$$(2) \qquad r_n(\varepsilon) \leq l = \sqrt{d - 1 + \left(1 - \frac{(\log \varepsilon^{-1})^{1/d}}{n^{1/d}}\right)^2}.$$

Putting (1) and (2) together:

$$\delta(n, \varepsilon) = r_n(1-\varepsilon) - r_n(\varepsilon) \geq u - l$$

$$\geq \sqrt{d}\left(1 - \frac{c\varepsilon^{1/d}}{n^{1/d}}\right) - \sqrt{d - 1 + \left(1 - \left(\frac{\log \varepsilon^{-1}}{n}\right)^{1/d}\right)^2}$$

$$= \sqrt{d}\left(1 - \frac{c\varepsilon^{1/d}}{n^{1/d}}\right) - \sqrt{d - 2\left(\frac{\log \varepsilon^{-1}}{n}\right)^{1/d} + \left(\frac{\log \varepsilon^{-1}}{n}\right)^{2/d}}$$

$$= \sqrt{d}\left[\left(1 - \frac{c\varepsilon^{1/d}}{n^{1/d}}\right) - \sqrt{1 - \frac{2}{d}\left(\frac{\log \varepsilon^{-1}}{n}\right)^{1/d} + \frac{1}{d}\left(\frac{\log \varepsilon^{-1}}{n}\right)^{2/d}}\right]$$

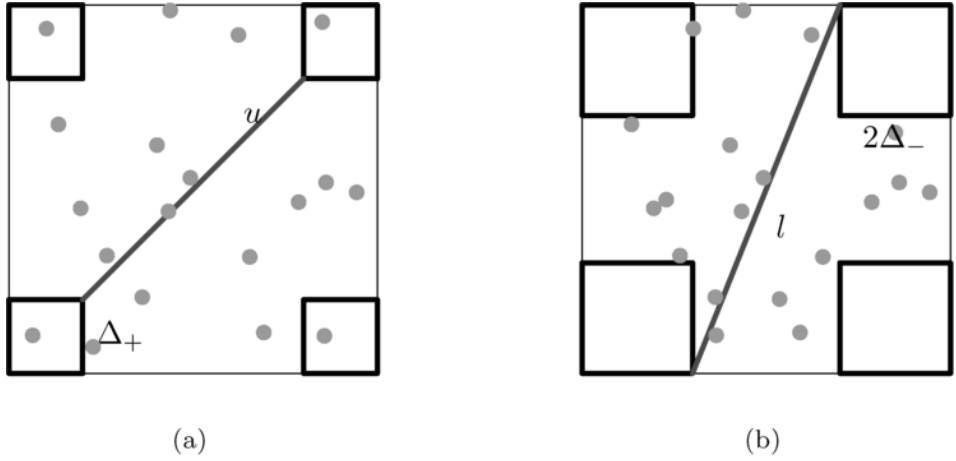

FIG. 1. *Definition of $\Delta_+$ and $\Delta_-$.*

$$= \sqrt{d}\left[\left(1 - c\left(\frac{\varepsilon}{n}\right)^{1/d}\right) - 1\right.$$
$$+ \frac{1}{2}\frac{2}{d}\left(\frac{\log \varepsilon^{-1}}{n}\right)^{1/d}\left(1 - \frac{1}{2}\left(\frac{\log \varepsilon^{-1}}{n}\right)^{1/d}\right)$$
$$\left. + \frac{1}{8}\frac{4}{d^2}\left(\frac{\log \varepsilon^{-1}}{n}\right)^{2/d}\left(1 - \frac{1}{2}\left(\frac{\log \varepsilon^{-1}}{n}\right)^{1/d}\right)^2 + \cdots\right]$$
$$= \sqrt{d}\left[\frac{\log^{1/d}\varepsilon^{-1} - cd\varepsilon^{1/d}}{dcn^{1/d}} + o\left(\frac{1}{n^{1/d}}\right)\right]$$
$$= \Omega\left(\frac{\log^{1/d}\varepsilon^{-1}}{n^{1/d}}\right). \qquad \square$$

Observe that for any $\kappa > 0$ constant, if $\varepsilon = n^{-\kappa}$, our lower bound for $d \geq 3$ matches our upper bound on the threshold width, and is only a factor of $\Theta(\log^{1/4} n)$ away for $d = 2$. For any constant $\varepsilon$, the difference between the bounds is $\mathcal{O}(\log^{1/d} n)$ and $\mathcal{O}(\log^{3/4} n)$, respectively.

**4. Bounding the bottleneck matching for $d \geq 3$.** We now present a simpler proof of $M_n = \Theta(r_c)$ when $d \geq 3$ than the proof presented in [26]. We emphasize that even though we also recursively subdivide the cube, our principle is different. In our proof, at the end of each step, the points are distributed uniformly in each subcuboid. This requires careful choice of the cutting plane and a linear transformation based on order statistics. This, however, permits us to bound the matching length via Chernoff bounds, as opposed to estimating the aspect ratios of rectangular solids as in [26].

The basic idea is to divide the unit square into $n$ equal boxes. Given $n$ points distributed uniformly on the square, move the points so that there is roughly a single point in each box. Now consider the two samples of red and blue points. Apply this process to both samples. Let $\Delta$ be the maximum distance by which a red point is shifted, and similarly, let $\Delta'$ be the maximum shift for any blue point, along any coordinate direction. Then the triangle inequality tells us that the bottleneck matching is less than $\sqrt{d}(\Delta + \Delta' + n^{-1/d})$.

We shall now use this idea to bound the bottleneck matching in $[0,1]^d$.

To move each point into its unique box we follow a recursive process. We shall provide only an informal description here, relegating the details to the Appendix. Moreover, for simplicity, we shall suppose $n$ to be a power of 2. First divide the square by drawing a vertical line so that there are exactly $n/2$ points in each half. Transform the $x$-coordinates of the points in each half so that they are uniformly distributed in $[0, 1/2]$ and $[1/2, 1]$. Now repeat the process along the $y$-axis for each rectangle, and then along the $z$-axis and so on. Repeat when all coordinate axes have been done once, and so on. One can carry this process for $\log n$ steps so that there is exactly one point in each box. However, for $d \geq 3$ it is better to stop at the $j$th step, where $j < \log n$. Choose $j = \lceil d^{-1} \log_2(n/\log n) \rceil$. Then the side of the box and hence the shift thereafter is $\leq 2^{-j}$. With this observation we can now show

PROPOSITION 4.1. *If $M_n$ is the weight of the bottleneck matching on a geometric random bipartite graph on $2n$ points in $[0,1]^d$, for $d \geq 3$, then for any $\beta > 1$, we can find a constant $c_d > 0$ such that*

$$\mathbb{P}\{M_n > c_d \beta r_c\} \leq \frac{1}{n^{\beta^2 - 1}},$$

*so that $M_n = \mathcal{O}(r_c)$ w.h.p.*

PROOF. To estimate $M_n$, we compute the total shift experienced by an arbitrary point. To this end, we shall find the shift along each axis, and so shall concentrate on one coordinate at a time. Let $x_1, \ldots, x_d$ denote the coordinates. We regard a *step* as one cycle in which divisions along all axes have been completed, according to the scheme described above. Therefore, if a step begins with a $d$-dimensional cube of side $l$ containing $n$ points, by the end of the step, the cube has been divided into $2^d$ cubes of side $l/2$, with $n/2^d$ points each.

Let $n_i = n/2^{d(i-1)}$ denote the number of points in a subcube at the beginning of the $i$th step, and let $l_i = 2^{-i+1}$ denote the length of the side of the cube. Lemma A.2 implies that for any point in the left half of such a subcube, the shift $\delta_i^{(k)}$ in the $x_k$-direction experienced during the $i$th step satisfies

$$\mathbb{P}\{|\delta_i^{(k)}| > \gamma_i\} \leq 2\exp\left(-\frac{\gamma_i^2}{l_i^2}n_i\right) \qquad \text{for any } \gamma_i > 0.$$

Therefore, if $\delta_i$ is the total shift suffered by a point during the $i$th step

$$\mathbb{P}\{|\delta_i| > \gamma_i\} \leq \mathbb{P}\left\{\max_{1\leq k\leq d}|\delta_i^{(k)}| > \gamma_i\right\} \leq \sum_{k=1}^d \mathbb{P}\{|\delta_i^{(k)}| > \gamma_i\}$$

$$\leq 2d\exp\left(-\gamma_i^2 \frac{n}{2^{(i-1)(d-2)}}\right).$$

Now fix a $\beta > 1$, and choose $\gamma_i$ such that $\gamma_i^2 \cdot n/2^{(i-1)(d-2)} = \beta^2 \log n$. Observe that with this choice, $\gamma_i$ is decreasing with $i$ only when $d > 2$. Let $\Delta$ be the maximum total shift experienced by any point. Then it must be that

$$\mathbb{P}\left\{|\Delta| > \sqrt{d}\sum_{i=1}^j \gamma_i\right\} \leq 2dnj\exp(-\beta^2\log n),$$

which follows from taking the union bound over all $n$ points, and all $d$ coordinates, and the fact that after $j$ steps, we divided a given coordinate at most $j$ times. Notice that after $j = \log n$ steps, the side of the subcube reduces to $2^{-\log n+1} = \Theta(1/n)$, and therefore subsequent shifts cannot move the point by more than $\mathcal{O}(1/n)$. Hence we can halt the subdivisions after $\log n$ steps, knowing that the matched point is already within $\mathcal{O}(1/n)$. Hence,

$$\mathbb{P}\left\{|\Delta| > \sqrt{d}\sum_{i=1}^{\log n} \gamma_i\right\} \leq 2dn\log n\exp(-\beta^2\log n).$$

However, with a little bit of work, one can see that

$$\sum_{i=1}^{\log n}\gamma_i = \beta\sqrt{\frac{\log n}{n}}\sum_{i=1}^{\log n} 2^{(i-1)(d-2)/2}$$

$$\leq 2\cdot\beta\sqrt{\frac{\log n}{n}}2^{(d-2)/2(\log_2 n - \log_2\log n)/d}$$

$$= 2\beta\left(\frac{\log n}{n}\right)^{1/d}.$$

Recall that $(\log n/n)^{1/d} = r_c \pi_d^{1/d}$, so that we have just shown

$$\sum_{i=1}^j \gamma_i \leq 2\beta\pi_d^{1/d}r_c.$$

Therefore, we have

$$\mathbb{P}\{|\Delta| > 2\sqrt{d}\beta\pi_d^{1/d}r_c\} \leq \mathbb{P}\left\{|\Delta| > \sqrt{d}\sum_{i=1}^{j}\gamma_i\right\} \leq \frac{2d\log n}{n^{\beta^2-1}}.$$

After $j$ steps the side of the cube is $2^{-j}$ and hence if we arbitrarily move points within the subcube, the extra shift is at most $\sqrt{d}2^{-j}$. Therefore, if we choose $c_d$ to be any constant larger than $\sqrt{d} + 2\sqrt{d}\pi_d^{1/d}$, we get $|\Delta| \leq c_d\beta r_c$ with probability at least $1 - n^{1-\beta^2}$. $\square$

We note in passing that for the above method to provide a bound in the $d=2$ case, one must proceed for $\log n$ steps, so that there is only a single point in each box. However, in this case, one only gets an $\mathcal{O}(r_c \log n)$ bound, which is off by $\log^{1/4} n$ from the sharp bound in [19].

**5. The bound for $d=1$.** Given $n$ points uniformly distributed in $[0,1]$, follow the recursive division procedure described in the last section. In this case, at each step the number of points decreases by half. Therefore, we obtain a stronger result:

PROPOSITION 5.1. *For any $\beta > 0$*

$$\mathbb{P}\{M_n^{(1)} \geq \beta/\sqrt{n}\} = \mathcal{O}(\exp(-c\beta^2))$$

*for some positive constant $c$.*

PROOF. For the sake of simplicity we assume that $n = 2^k$ for some $k \in \mathbb{N}$; this makes no difference to the proof except for simplifying some expressions. In the $i$th step there are $2^i$ sets of $n/2^i$ points each. Therefore, if $\delta_i$ is the shift of an arbitrary set, then for $\beta_i > 0$, by Lemma A.2:

$$\mathbb{P}\left\{|\delta_i| \geq \frac{\beta_i}{\sqrt{n}}\right\} \leq 2\exp(-2^i\beta_i^2),$$

so that the maximum shift $\Delta$ of any point satisfies

$$\mathbb{P}\left\{\max|\Delta| \geq \frac{\sum_i \beta_i}{\sqrt{n}}\right\} \leq 2\sum \frac{n}{2^i}\exp(-2^i\beta_i^2).$$

Choose the $\beta_i$'s such that $2^i\beta_i^2 = \beta_0^2 + i$, for some $\beta_0$. Then $\sum_i \beta_i \leq K\beta_0$ for some suitable constant $K > 0$. Taking $\beta = K\beta_0$, we get

$$\mathbb{P}\left\{M_n^{(1)} \geq \frac{\beta}{\sqrt{n}}\right\} \leq c'\exp(-c\beta^2),$$

for some constants $c, c' > 0$. $\square$

**6. Conclusion.** We have shown that all monotone graph properties have a sharp threshold for random geometric graphs. Moreover, this threshold is sharper than the one for Bernoulli random graphs.

We have a sharp result for $d = 1$. For the $d \geq 3$ we believe the upper bound of $\mathcal{O}(r_c)$ to be actually tight. For the $d = 2$ case we believe the upper bound to be $\mathcal{O}(r_c)$ as well, though this cannot be obtained via bottleneck matchings.

## APPENDIX: ESTIMATING THE SHIFT IN EACH RECURSIVE STEP

To establish a bound on the amount by which each point is moved, we must examine the "shrinking" and "stretching" process formally. For simplicity we concentrate on the $d = 2$ case. Ignore the $y$-coordinates. Then we have $n$ i.i.d. $\mathsf{Unif}[0,1]$ points $x_1, \ldots, x_n$. It is well known [22] that the $k$ smallest points are i.i.d. uniform in $[0, x^{(k+1)})$, and that the $n - k$ largest points are i.i.d. uniform in $(x^{(k)}, 1]$, where $x^{(i)}$ is the $i$th order statistic. Suppose that $n$ is even—the analysis below applies *mutatis mutandis* when $n$ is odd. Set $\delta_l = x^{(n/2+1)} - 1/2$ and $\delta_r = 1/2 - x^{(n/2)}$. Then transform the points as follows:

$$x^{(i)} \mapsto \begin{cases} x^{(i)} \dfrac{1/2}{1/2 + \delta_l}, & \text{for } i = 1, \ldots, \dfrac{n}{2}, \\ 1 - (1 - x^{(i)}) \dfrac{1/2}{1/2 + \delta_r}, & \text{for } i = \dfrac{n}{2} + 1, \ldots, n. \end{cases}$$

This transform leaves the smallest $n/2$ points uniformly distributed in $[0, 1/2]$ and the largest $n/2$ points uniformly distributed in $[1/2, 1]$. This process is now repeated $\lceil \log n \rceil$ times alternating the $x$- and $y$-coordinates. The maximum shift at any step is not more than $|\delta_l|$ for the points on the left and not more than $|\delta_r|$ for points on the right. We shall use $\delta = \max\{\delta_l, \delta_r\}$.

Let $X_t$ be the number of points in $[0, t]$ prior to the transformation. The following result is immediate:

LEMMA A.1. *For $0 < \gamma < 1/2$,*

$$\mathbb{P}\{|\delta| > \gamma\} \leq 2\mathbb{P}\left\{X_{1/2+\gamma} < \frac{n}{2}\right\}.$$

To bound the last probability, observe that $X_t$ is just the sum of $n$ i.i.d. Bernoulli's that are 1 with probability $t$. Hereafter $\beta > 0$ is some constant.

LEMMA A.2. *The shift $\delta$ of any point in the recursion step satisfies*

$$\mathbb{P}\{|\delta| > \gamma\} = \mathcal{O}(\exp(-n'(\gamma/l)^2)),$$

*where $n'$ is the number of points in the subcuboid being divided, and $l$ is length of the side that is being divided.*

PROOF. The proof is a straightforward application of Chernoff's bound. Assume wlog that $l = 1$:

$$\begin{aligned}
\mathbb{P}\{|\delta| > \gamma\} &\leq 2\mathbb{P}\bigg\{X_{1/2+\gamma} < \frac{n'}{2}\bigg\} \\
&= 2\mathbb{P}\bigg\{X_{1/2+\gamma} < n'\bigg(\frac{1}{2} + \gamma\bigg) - n'\gamma\bigg\} \\
&\leq 2\exp\bigg(-\frac{n'^2\gamma^2}{2n'(\gamma + 1/2)}\bigg) \\
&\leq 2\exp(-n'\gamma^2) \qquad \text{since } \gamma \leq 1/2. \qquad \square
\end{aligned}$$

Generalization to the $d \geq 3$ case is straightforward.

**Acknowledgments.** We should like to thank Peter Glynn, David Kempe, Rajeev Motwani and Devavrat Shah for useful discussions. We are especially grateful to Greg McColm for his comments on an earlier draft of the paper.

A. GOEL
DEPARTMENT OF MANAGEMENT SCIENCE
  AND ENGINEERING
AND (BY COURTESY)
DEPARTMENT OF COMPUTER SCIENCE
STANFORD UNIVERSITY
STANFORD, CALIFORNIA 94305
USA
E-MAIL: ashishg@stanford.edu

S. RAI
DEPARTMENT OF MANAGEMENT SCIENCE
  AND ENGINEERING
STANFORD UNIVERSITY
STANFORD, CALIFORNIA 94305
USA
E-MAIL: sanat@stanford.edu

B. KRISHNAMACHARI
DEPARTMENT OF ELECTRICAL ENGINEERING
  AND DEPARTMENT OF COMPUTER SCIENCE
UNIVERSITY OF SOUTHERN CALIFORNIA
LOS ANGELES, CALIFORNIA 90089
USA
E-MAIL: bkrishna@usc.edu